%Paper: alg-geom/9601016
%From: Tohru Nakashima <nakasima@math.metro-u.ac.jp>
%Date: Thu, 18 Jan 96 00:25:46 JST

\documentstyle{amsppt}

\pageheight{17cm}
\pagewidth{13.3cm}
\vcorrection{0.8cm}
\magnification=\magstep1
\topmatter
\title
Space of conformal blocks in 4D WZW theory
\endtitle

\author
Tohru Nakashima \\
Yuichiro Takeda
\endauthor
\address
Department of Mathematics,
Tokyo Metropolitan University,
Minami-Ohsawa 1-1, Hachioji-shi,
Tokyo, 192-03 Japan
\endaddress
\email
nakasima\@math.metro-u.ac.jp,
takeda\@math.metro-u.ac.jp
\endemail
\keywords
conformal blocks, stable bundles, moduli spaces
\endkeywords
\subjclass
14D20, 14J60
\endsubjclass
\abstract
An algebro-geometric definition of the space of conformal blocks in the 4D WZW
theory is proposed. As an example, the dimensions of these spaces are
calculated for ruled surfaces.
\endabstract

\endtopmatter

\redefine\red{\operatorname{red}}
\redefine\Quot{\operatorname{Quot}}
\redefine\dim{\operatorname{dim}}

\redefine\PGL{\operatorname{PGL}}

\redefine\Det{\operatorname{Det}}

\redefine\max{\operatorname{max}}

\redefine\Tr{\operatorname{Tr}}
\redefine\SL{\operatorname{SL}}

\subheading{\bf 1. Introduction }

\vskip 1pc

Two-dimensional conformal field theory has attracted attention of many
mathematicians in the last years.
Especially the Verlinde conjecture giving the dimensions of the spaces of
conformal blocks has been at the center of interest for its remarkable
connection with various fields of mathematics including algebraic geometry and
representation theory.
At present this conjecture has been proved rigorously [1]. \par
Recently, a four dimensional analogue of Wess-Zumino-Witten theory has been
proposed [3].
As in 2D CFT, one of the basic objects of study in this theory is the space of
conformal blocks, which is described as the space of sections of a certain line
bundle on the moduli of holomorphic vector bundles on an algebraic surface.
In loc.cit.an analogue of the Verlinde conjecture for the dimensions of these
spaces has also been presented.
Mathematically, however, this conjecture is not well defined since the relevant
moduli space is rarely compact.
The purpose of this note is to interpret the space of 4D conformal blocks by
means of natural line bundles on the moduli of Gieseker semistable torsion-free
sheaves, which have been defined by one of the present authors [4].
As an example, we calculate the dimension of the space of global sections of
these line bundles for rational ruled surfaces.

\vskip 1pc

\subheading {\bf 2. The space of conformal blocks }

\vskip 1pc

Let $X$ be a smooth projective surface over $\Bbb C$ and $H$ an ample line
bundle on $X$.
We denote by $\Cal M=\Cal M_H(r,c_1,c_2)$ the moduli space of rank $r$ vector
bundles with Chern classes $c_1$,\,$c_2$ which are $H$-stable in the sense of
Mumford-Takemoto.
Although this is not in general a projective variety, it has a compactification
$\Cal M^{ss}$, the moduli of Gieseker semistable torsion-free sheaves [2].
Let $\overline {\Cal M}$ denote the closure of $\Cal M$ in $\Cal M^{ss}$ for
Zariski topology, which can be regarded as a minimal compactification of $\Cal
M$.
We recall below the construction of determinant line bundles on $\overline
{\Cal M}$ as defined in [4]. \par
We fix a sufficiently large integer $N$ and let $\Cal O_X(N)=H^{\otimes N}$.
For a suitable integer $P(N)$ depending on $N$,
$\Cal M^{ss}$ is contructed as the geometric invariant theory quotient of a
Quot scheme $\Cal Q=\Quot ^{ss}(\Cal O_X(-N)\otimes \Bbb C^{P(N)})$ under $\SL
(P(N))$ action. We recall that $\Cal Q$ parametrizes the quotient sheaves $E$
$$
\Cal O_X(-N)\otimes \Bbb C^{P(N)}\to E
$$
together with an isomorphism $\Bbb C^{P(N)}\cong H^0(X,E(N))$. \par
We fix a universal family $\Cal F$ on $X\times \Cal Q$. If $C$ is a smooth
curve on $X$ such that the intersection number $c_1\cdot C$ is divisible by
$r$, we choose a line bundle $L$ on $C$ satsfying $\chi (\Cal F_x^C\otimes
L)=0$ for every sheaf $\Cal F_x$ corresponding to a point $x\in \Cal Q$.
We define a line bundle $\Det _{\Cal F}(C)$ on $\Cal Q$ as follows.
$$
 \Det _{\Cal F}(C)=\det ((p_{\Cal Q})_!(\Cal F^C\otimes p_C^*L))^{\vee }
$$
where $\Cal F^C=\Cal F_{\vert C\times \Cal Q}$ and $p_C$, $p_{\Cal Q}$ denote
the projections from $C\times \Cal Q$ to $C$ and $\Cal Q$, respectively.
Since $\PGL(P(N))$ acts trivially on the fiber of $\Det _{\Cal F}(C)$ at every
point with closed orbit, we can descend the line bundle to $\Cal M^{ss}$ which
will be also denoted by $\Det _{\Cal F}(C)$.
 It can be easily seen that this construction does not depend on the choice of
$L$.
Let $\Cal L_{\Cal F}(C)=\Det _{\Cal F}(C)_{\vert \overline {\Cal M}}$ be the
restricted line bundle. We define {\it the space of conformal blocks} as
follows
$$
  Z_{\Cal F,C}(r,c_1,c_2)=H^0(\overline {\Cal M},\Cal L_{\Cal F}(C)).
$$
If there exists a universal family $\Cal E$ on $X\times \overline {\Cal M}$, we
set
$$
 \Cal L_{\Cal E}(C)=\det ((p_{\overline {\Cal M}})_!(\Cal E^C\otimes
p_C^*L))^{\vee }
$$
and similarly we define the space $Z_{\Cal E,C}(r,c_1,c_2)$. \par
To see the relation of these spaces with 4D WZW theory, we assume that $X$ is a
regular algebraic surface( i.e. $H^1(X,\Cal O_X)=0$ ), $c_1=0$ and that there
exists a universal family $\Cal E$ on $X\times \overline {\Cal M}$.
Then the Grothendieck-Riemann-Roch theorem yields
$$
     c_1(\Cal L_{\Cal E}(C))=c_2(\Cal E)/[C]
$$
where $[C]\in H_2(X,\Bbb Z)$ denotes the homology class of $C$ and $/$ is the
slant product. This may be considered as a higher rank generalization of
Donaldson's $\mu $ map. \par
We denote by $\overline {\Cal M}_{\red }$ the scheme theoretic reduction of
$\overline {\Cal M}$. Let $F$ be the curvature form of a hermitian connection
on $\Cal E$ compatible with the complex structure.
If $[C]$ is the Poincar\'e dual of a holomorphic two form $\omega $, then
$c_1(\Cal L_{\Cal E}(C))$ is represented on the smooth locus of $\overline
{\Cal M}_{\red }$ by the following two form [5]
$$
   \frac {1}{8\pi ^2}\int _X\Tr (F\wedge F)\wedge \omega
$$
which essentially reproduces the line bundle considered by physicists in [3].
This suggests that $Z_{\Cal F,C}(r,c_1,c_2)$ should be a mathematical
formulation of the space of conformal blocks in 4D WZW theory.
An interesting problem of computing its dimension remains open at present,
although for bundles with $c_1=0$ a formula of Verlinde type has been presented
[3]. Inspired by the conjecture, we shall compute in the next section the
dimension of $Z_{\Cal F,C}(2,c_1,c_2)$ for certain bundles with $c_1\ne 0$ on a
rational ruled surface.

\vskip 1pc

\subheading{\bf 3. An example }

\vskip 1pc

 Let $\pi :X\to \Bbb P^1$ be a ruled surface, namely the $\Bbb P^1$-bundle
associated to some rank two vector bundle on $\Bbb P^1$.
We denote by $\Sigma $, $f$ the divisor class of a section and a fiber
respectively.
Let $e=-\Sigma ^2$ and assume that $c>\max (-e/4,0)$.
Then by [6], we can find an ample line bundle $H_c$ such that a rank two bundle
$E$ with $c_1(E)=\Sigma $, $c_2(E)=c$ is $H_c$-stable if and only if $E$ is
given as the following nontrivial extension of line bundles
$$
 0\to \Cal O(\Sigma -cf)\to E\to \Cal O(cf)\to 0.
$$
Thus the moduli space $\Cal M=\Cal M_{H_c}(2,\Sigma ,c)$ is particularly simple
since $H$-stability and Gieseker semistability coincide and torsion-free
sheaves do not appear: we have $\Cal M=\overline {\Cal M}=\Cal M^{ss}$. Since
the extensions as above are parametrized by the vector space of dimension
$4c+e-2$
$$
V_c=H^1(X,\Cal O(\Sigma -2cf)),
$$
we conclude that $\Cal M$ is isomorphic to the projective space $\Bbb
P(V_c^{\vee })\cong \Bbb P^n$ where $n=4c+e-3$.
Furthermore, it is easy to see that there exists a universal rank two bundle
$\Cal E$ on $X\times \Cal M$ which is given as the following extension
$$
 0\to p_X^*\Cal O(\Sigma -cf)\to \Cal E\to p_X^*\Cal O(cf)\otimes p_{\Cal
M}^*\Cal O(1)^{\vee }\to 0
$$
where $p_X$, $p_{\Cal M}$ denote the projections from $X\times \Cal M$ to $X$,
$\Cal M$ and $\Cal O(1)$ is the tautological line bundle on $\Cal M=\Bbb P^n$.
\par
Let $C$ be a smooth curve in $X$ such that $C\cdot \Sigma $ is even.
Using the universal extension above, a straightforward calculation shows that
$\Cal L_{\Cal E}(C)$ is isomorphic to $\Cal O(m)$ where $m$ is the following
intersection number
$$
  m=\left (cf-\frac {\Sigma }{2}\right )\cdot C.
$$
Therefore we obtain
$$
\dim Z_{\Cal E,C}(2,\Sigma ,c)=\dim H^0(\Bbb P^n,\Cal O(m))=\binom {n+m}{n}.
$$

\enddocument